\documentclass[12pt]{amsart}
\usepackage{amssymb,amsfonts,amsmath}
\usepackage{times}
\usepackage{graphicx}

\begin{document}

\title[DBT\ \&\ CA]
{Discrete Baker Transformation and Cellular Automata}

\author{Valeriy K. Bulitko}
\email{booly@shaw.ca}

\vspace*{-2.0cm} \sloppy \thispagestyle{empty}

\maketitle

%\ \footnote{This work was partially supported by the
%research grant "Families of Polynomials Associated to
%Certain Discrete Dynamical Systems" of Athabasca
%University.}

\begin{abstract}
In this paper we propose a rule-independent description of
applications of cellular automata rules for one-dimensional additive cellular
automata on cylinders of finite sizes.
This description is shown to be a useful tool for for answering questions
about automata's state transition diagrams (STD). The approach
is based on two transformations: one (called {\sl
Baker transformation}) acts on the $n$-dimensional
Boolean cube $\frak B^n$ and the other (called {\sl
index-baker transformation}) acts on the cyclic group of
power $n$. The single diagram of Baker transformation
in $\frak B^n$ contains an important information about
all automata on the cylinder of size $n$. 
Some of the results yielded by this approach can be viewed as a generalization and extension of certain results by O. Martin, A. Odlyzko, S. Wolfram~[1]. 
Additionally, our approach leads to a convenient 
language for
formulating properties, such as possession of cycles
with certain lengths and given diagram heights,
of automaton rules. 
\end{abstract}

{\footnotesize \tableofcontents}

\newtheorem{theo}{Theorem}
\newtheorem{lem}{Lemma}
\newtheorem{cor}{Corollary}
\newtheorem{prop}{Proposition}
\newtheorem{ex}{Example}
\newtheorem{rem}{Remark}
\newtheorem{rul}{Rule}
\newtheorem{conj}{Conjecture}

\section{Introduction}

In this paper we consider the behavior of one-dimensional
cellular automata acting on a finite cylinder of size
$n$. The idea of our approach is as follows.
Maximum and minimum of descriptive complexity of the rules produced by iterative applications of an arbitrary rule depend, in principle, on the cylinder size and the generating rule. It turns out, however, that in the case of additive finite one-dimensional automata, the minimum complexity is reached at the second iteration regardless of the cylinder size and the generating rule. We show this by introducing a straight-forward rule-independent procedure that yields the results of the second iteration without applying the generating rule twice. This procedure is called {\sl Discrete Baker Transformation} (DBT) and found to be a useful tool for answering a broad spectrum of questions on automata behavior.

In the following we will introduce the notation used in the rest of the paper. A cellular automaton is called additive if it is defined by additive rules ($X$) acting on the cell and its right neighbors [1, 2].
In this paper we consider one-dimensional
additive cellular automata (ACA) on a finite cylinder of
size $n$ with the states from boolean cube $\frak B^n$.
Let us call an automaton state {\em cyclic} if it belongs to
some cycle of the state diagram of the automaton.
We denote the length of a string $w$ by $|w|$ and
the parity of a binary string $w$ by $\varrho(w)$.
Pair $(n,X)$ where $n\in\mathbb N$ and $X$ is a
finite 0,1-sequence of length $m, m\le n$, defines an
additive cellular automaton (ACA) on the cylinder of
size $n$. We denote the automaton by $\frak A(n,X)$. The {\em standard
length} of rule $X$ for an automaton on the cylinder of
size $n$ is $n$. Thus, the short notation 
$X, |X|=m<n$ for a automaton rule on an 
$n$-cylinder means that $X$ must be padded with zeroes
at its right end until the length of $n$. The rest of the paper uses the full notation by default.

The states of $\frak A(n,X)$ constitute a boolean cube $\frak B^n$ of dimension $n$. We
write the strings as words in the alphabet $\{0,1\}$
or as vectors (e.g., $[i_1,i_2,\dots,i_n]$).
For any given $n$ two particular strings play an
important role. They are $0^n\ (=[0,\dots,0]),1^n\
(=[1,\dots,1])$. We refer to them as $\mathbf 0,\mathbf 1$ correspondingly.
Rule $X$ acts as a linear operator in space of all
states $\frak B^n$ of $\frak A(n,X)$. We denote the
operator by $\hat X$. Accordingly, $\hat X*s$ is the result of
application of the rule $X$ to state $s$. Rules themselves (in
the standard form) comprise the boolean cube $\frak B^n$.

Any non-cyclic ACA state evolves into a cyclic state in finite time. The maximum time of this evolution is called the {\em height} of the ACA and is denoted by $h^{\star}$. It is defined to be 0 if the ACA has cyclic states only. To specify a concrete 
automaton $\frak A(n,X)$ we use $h^{\star}(n,X)$.
Another important attribute of an ACA is the
distribution of cycle lengths and, in particular,
the maximum cycle length.

\section{Iterations of $\hat X$}

Clearly, the characteristics of STD s.t. maximal
height or cycle length relate to lengths of the
oriented chains of the correspondent graph. These
chains of states in STD are produced by sequential
actions of the same operator $\hat X$, i.e. iterations
of the action or operators which are the degrees $\hat
X^i$ of given an operator $\hat X$. Because we denote
by $\hat I$ the identical operator and accept $\hat
X^0=\hat I$, the index of degree
$\hat X^i$ may run all natural numbers.\\

\subsection{Kernels and images of $\hat X$}
\  \\
The following well known for all linear spaces and
linear operators in them\footnote{Actually - for all
commutative groups and their homomorphisms.} facts we
apply below for operator $\hat X$ and space $\frak
B^n$.

\begin{prop}\label{prop_ab_kernels}

The following statements are well known:
\begin{enumerate}

\item The intersection of subspaces $\frak B^n$ is a
subspace, and $\{\mathbf 0\}$ is the least subspace of
$\frak B^n$ with respect to inclusion.

\item Image $\mathrm{Im}(\hat X)$ and kernel
$\ker(\hat X)$ are linear subspaces of $\frak B^n$.

\item $\forall i[\mathrm{Im}(\hat
X^{i+1})\subseteq\mathrm{Im}(\hat X^i)\ \&\ \ker(\hat
X^{i+1})\supseteq\ker(\hat X^i)]$.

\item $\forall i[\mathrm{Im}(\hat
X^{i+1})\subsetneq\mathrm{Im}(\hat X^i)\iff\ker(\hat
X^{i+1})\supsetneq\ker(\hat X^i)\iff\mathrm{Im}(\hat
X^i)\cap\ker(\hat X^i)\neq\{\mathbf 0\}]$.

\item $\forall i[\ker(\hat X^{i+1})=\ker(\hat
X^i)\implies\ker(\hat X^{i+2})=\ker(\hat X^{i+1})]$.

\item $\forall i[|\mathrm{Im}(\hat
X^{i+1})|=\frac{|\mathrm{Im}(\hat X^{i})|}{|\ker(\hat
X^i)\cap\mathrm{Im}(\hat X^{i})|}]$.

\end{enumerate}
\end{prop}
 \  \\

\subsection{Composition of circulants}
\ \\

Let $L,L=[a_0,\dots,a_{n-1}],$ is $0,1$-sequence of
length $n$ and $\sigma$ is the cyclic shift right of
$L$, i.e. $\sigma*L=[a_{n-1},a_0,\dots,a_{n-2}]$. Also
let $s^{\circlearrowleft}$ denote string $s$ being
read from the end to the beginning (i.e. in reverse
order).

%$\stackrel{\curvearrowleft}{s}$
%$\stackrel{\leftarrow}{s}$

We call the square matrix which $i$-th row is
$\sigma^{i-1}(L)$ by circulant of $L$ and denote as
$C(L)$. Since circulant is completely defined by its
the first row, it has a sense to call the first row of
a circulant {\sl leader row} or {\sl leader}. Let us
define for any $n$ a binary operation $\boxtimes:\frak
B^n\times\frak B^n\to\frak B^n$ on leaders of
circulants as following:
\begin{equation}~\label{pre-bak}
(L\boxtimes M)_j=\sum_{k=0}^{n-1}L_k
 (M^{\circlearrowleft})_{k-j-1(\mathrm{mod}\ n)},
 j=\overline{0,n-1}.
\end{equation}
It can be noted that
$(M^{\circlearrowleft})_{k-j-1(\mathrm{mod}\
n)}=M_{j-k(\mathrm{mod}\ n)}$. We use expression
$M^{\circlearrowleft}$ to present a "geometric"
structure of the operation: the second operand $M$ in
first becomes reversed and then cyclic shifted on one
position right.

The next lemma shows the meaning the operation we
introduced by~(\ref{pre-bak}).

\begin{lem}\label{lem.3}
Let $|L|=|M|=n$. Then $C(L)C(M)=C(Q)$ where $Q=
L\boxtimes M$.
\end{lem}

{\bf Proof.} Let $L^i$ be the $i-th$ row of $C(L)$ and
$T^j$ be the $j$-th column of $C(M)$.

Since $L^i=\sigma^i*L$ and
$T^j=\sigma^{j+1}*M^{\circlearrowleft}$, then
\begin{equation*}
(C(L)C(M))_{i,j}=\sum_{k=0}^{n-1}(\sigma^i*L)_k
(\sigma^{j+1}*M^{\circlearrowleft})_k=\sum_{k=0}^{n-1}
L_{k-i(\mathrm{mod}\
n)}(M^{\circlearrowleft})_{k-j-1(\mathrm{mod}\ n)}.
\end{equation*}
So $Q$ is the first ($i=0$) row of the matrix product.
From here we get
$$
(C(L)C(M))_{i+1,j+1}=\sum_{k=0}^{n-1}
L_{k-i-1(\mathrm{mod}\
n)}M^{\circlearrowleft}_{k-j-2(\mathrm{mod}\ n)}.
$$
Setting $k'=k-1(mod\ n)$ and taking in account that
$k,k'$ are tied variables running the same scope of
numbers, we conclude that
$$
(C(L)C(M))_{i,j}=(C(L)C(M))_{i+1,j+1}.
$$
The last means that $C(L)C(M)$ is the circulant of
$Q$. $\ \ \ \Box$

Thus, $\boxtimes$ is just an image of matrix
multiplication of circulants in space of their
leaders. And because matrix multiplication is
associative, $\boxtimes$ is the associative operation
too.

This is an useful property of the operation:
\begin{lem}~\label{3a}
Suppose $|L|=|M|=qn'$ where $q,n'\in\mathbb N$. Then
\begin{eqnarray*}
\forall i[0\le i\le qn'-1\,\&\,q\nmid i\implies
L_i=M_i=0]\implies\\
\forall i[0\le i\le qn'-1\,\&\,q\nmid i\implies
(L\boxtimes M)_i=0].
\end{eqnarray*}
\end{lem}

{\bf Proof.} Because if $q\nmid k$ then $L_k=0$, we
can rewrite the definition~\ref{pre-bak} of
$\boxtimes$ in the next form:
\begin{equation}~\label{pre-bak+}
(L\boxtimes M)_j=\sum_{i=0}^{n'-1}L_{iq}
 M_{j-iq(\mathrm{mod}\ n)},
 j=\overline{0,n-1}.
\end{equation}
Now, clearly $q\nmid j\iff q\nmid (j-iq(\mathrm{mod}\
n))$ since $q|n$. From here and the condition of the
lemma for $M$ we draw that if $q\nmid j$ all
$M_{j-iq(\mathrm{mod}\ n)}$ are equal 0, and
$(L\boxtimes M)_j=0$. $\ \ \ \Box$

This technical lemma has an important meaning for the
theory. One application will be presented in section
about the reduction of the problem to compute
determinants (modulo 2) of automata rules.\\

\subsection{Discrete baker transformation (DBT)}
 \ \\

As we will see further, the case $L\boxtimes L$
("baker transformation") presents the special
interest.

In chaos theory the transformation $x_{n+1}=2\mu x_n$,
where $x$ is computed modulo 1 (for $\mu=1$ see ,
[4,p.272]) is called baker transformation or map. By
this analogy we will call {\sl discrete baker
transformation} on the set $\frak B^n$ of all boolean
$n$-tuples the mapping $\frak b:\frak B^n\to\frak B^n$
acting according to the rule:
\begin{equation}\label{baker_df}
(\frak b*L)_i=
\begin{cases}
\underset{j:\ 2j=i(\mathrm{mod}\ n)}{\bigoplus}L_j,&
\text{ if
}\{k|2k=i\ (\mathrm{mod}\ n)\}\neq\emptyset,\\
\ \ \ \ \ \ \ 0,&  \text{ else}.
\end{cases}
\end{equation}
\ \\
\begin{ex}\label{ex.5}
{\rm Let
$n=9,L=[a_0,a_1,a_2,a_3,a_4,a_5,a_6,a_7,a_8]$. Then
$\frak
b*L=%\linebreak[4]
[a_0,a_5,a_1,a_6,a_2,a_7,a_3,a_8,
a_4]$. In case $n=8$ we have $\frak
b*[a_0,a_1,a_2,a_3,a_4,a_5,%\linebreak[4]
a_6,a_7]=[a_0+a_4,
0, a_1+a_5, 0, a_2+a_6, 0, a_3+a_7, 0]$.}
\end{ex}

The next lemma proves the identity $L\boxtimes L=\frak
b*L$ and therefore explains the meaning of the
discrete baker transformation for us:

\begin{lem}\label{lem.4}
$C(L)C(L)=C(\frak b*L)$.
\end{lem}

{\bf Proof.} Let $|L|=n$. As we saw $C(L)C(L)=C(Q)$
where $Q_j=\sum_{k=0}^{n-1}L_k
 L^{\circlearrowleft}_{k-j-1(mod\ n)},
 j=\overline{0,n-1}.$
Because $L^{\circlearrowleft}_{k-j-1(mod\
n)}=L_{j-k(mod\ n)}$ we can write
$$
Q_j=\sum_{k=0}^{n-1}L_kL_{j-k(mod\ n)}.
$$
If we imagine string $L$ in the form of a ring then we
will see that $k$-th item in the sum is the product of
$k$-th component of $L$ counted from the $0$-th
component in positive direction (when numbers of
component increase) and $k$-th component counted from
$j$-th component in the opposite direction.

{\sl Case 1: $n$ is odd.} In this case for any $j$
there exists only one component equidistant from
$0$-th and $j$-th components. Its number $\gamma(j)$
is
\begin{center}
$\gamma(j)=
\begin{cases}
\frac{j}{2},& \text{ if $j$ is even,}\\
j+\frac{n-j}{2},&  \text{ if $j$ is odd}.
\end{cases}$
\end{center}
Since all these products of different components occur
in the sum twice, the sum (i.e. $Q_j$) is equal to
$L_{\gamma(j)}$. Since $2\gamma(j)=j(mod\ n)$ we have
$Q_j=(\frak b*L)_j,j=0,\dots,n-1$.

{\sl Case 2: $n$ is even.} This time either there
exist exactly two components
$L_{\gamma(j)},L_{\delta(j)}$ equidistant from $0$-th
and $j$-th components or no one at all. In the last
subcase $j$ is odd and the sum $Q_j$ equals to 0
because every item in the sum occur twice. In the
former subcase $j$ is even and
$\gamma(j)=\frac{j}{2},\delta(j)=\frac{n-j}{2}$. So we
have
\begin{center}
$Q_j=
\begin{cases}
0,& \text{ if $j$ is odd},\\
L_{\frac{j}{2}}+L_{\frac{n-j}{2}},&  \text{ if $j$ is
even}.
\end{cases}$
\end{center}
Again the result coincides with $\frak b*L$. $\ \ \
\Box$

The next result plays an important role in decoding of
the baker diagrams, see next sections.\\

\begin{theo}~\label{conserv_pr}
(i) {\bf Conservation principle:}
$\forall X[\det_2(|X|,X)=\det_2(|X|,\frak b*X)]$.\\
(ii) $\forall
X[\mathrm{rank}(C(X))\ge\mathrm{rank}(C(\frak b*X))]$.
\end{theo}

{\bf Proof.} (i) According to the lemma~\ref{lem.4}
$C(L)C(L)=C(\frak b*L)$. As it is well known
$\det(AB)=\det(A)\det(B)$. Therefore $\det(AB)$ is odd
iff the both of $\det(A),\det(B)$ are odd.

(ii) This statement is the prompt consequence of
$\frak b$ definition and
proposition~\ref{prop_ab_kernels}(3).
 $\ \ \ \Box$

Application a rule $\hat X$ to a state $s$ traces only
step in complete trajectory of the state. Two steps
produces the operator $(\hat X)^2$ or $\widehat{\frak
b*X}$; four step result is produced by $\widehat{\frak
b^2*X}$ and so forth, see the next theorem.

First, as usual we define
$$
\frak b^{i+1}*X=\frak b*(\frak b^i*X).
$$

\begin{theo} \label{th.5}
$\widehat{\frak b^i*X}=(\hat X)^{2^i}$.
\end{theo}

{\bf Proof.} For beginning we note that compositions
of operators of kind $\hat X$ is associative since
these are linear operators.

Induction on $i$. The basis $i=1$ was proved in
lemma~\ref{lem.4}. Now suppose the statement is true
for $i=k$ let's prove it for $i=k+1$. We have
$\widehat{\frak b^{i+1}*X}=_{\text{(by
definition)})}\widehat{\frak b*(\frak
b^i*X)}=_{\text{(by the induction
basis)}}(\widehat{\frak b^i*X})^2=_{\text{(by the
induction supposition)}}((\hat X)^{2^i})^2=_{\text{(by
associativity of the operator composition)}}\hat
X^{2^{i+1}}.\ \ \ \Box$

Despite the non-uniform scale, the result can play an
important role in studying asymptotic behavior of the
operators as $\hat X$ and therefore in understanding
the global structure of STD.

The amazing facts about behavior of $\frak b$ in
$\frak B^n$ is described by the next two lemmas. To
formulate and prove them it's worth to select "index
projection" of the baker transformation. For given a
number $n$ it is the map
$\natural:[0,\dots,n-1]\to[0,\dots,n-1]$ determined by
the rule
\begin{equation}\label{ind_baker_df}
\natural(i)=\mathrm{rem}(2i,n),
\end{equation}
(i.e. the remainder from division $2i$ by $n$). What
we mean talking about index projection is that $(\frak
b(L))_i=\underset{\natural(j)=i}{\oplus}L_j$,
see~(\ref{baker_df}). From here the following follows.

\begin{prop}\label{prop_ab_bak_and in-bak}
$\forall k\forall X
[\natural^k([0,\dots,|X|-1])=[0,\dots,|X|-1]\implies
\frak b^k(X)=X]$.
\end{prop}

{\bf Proof.} The condition
$\natural^k([0,\dots,|X|-1])=[0,\dots,|X|-1]$ implies
that $\natural$ is a permutation on $[0,\dots,|X|-1]$.
So $\frak b$ is a permutation of components of $X$
(nothing to glue). Therefore the implication is true;
yet the cycle of $\frak b$ could be even shorter:
0-vector is a good example. $\ \ \ \Box$

As usual one can draw diagrams of the mapping
$\natural$.

\begin{ex}\label{ex_ind_1}
{\rm The fig.~\ref{fig11} shows diagrams of $\natural$
on the segment $[0,\dots,32], n=33$. Note, $n$ is odd,
and the set of indices $\{0,\dots,32\}$ is partitioned
on 5 cycles: $\{0\}$, $\{11, 22\}$, $\{1, 2, 4, 8, 16,
32, 31, 29, 25, 17\}$, $\{3, 6, 9, 12, 15, 18, 21, 24,
27, 30\}$, $\{5, 7, 10, 13, 14, 19, 20, 23, 26, 28\}$
with lengths correspondingly 1,2,10,10,10.}
\end{ex}

\setlength{\unitlength}{.1cm}
\begin{figure}[htbp]
\begin{picture}(180,30)
\put(10,0){\includegraphics[width=10cm]{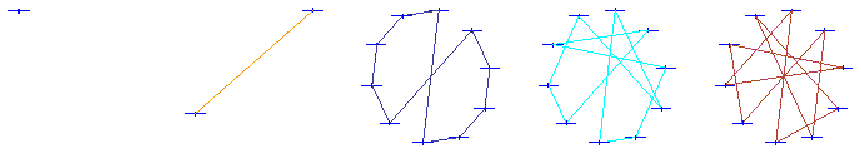}}
\end{picture}
\caption{The diagram of the mapping $\natural$ for
$n=33$.}\label{fig11}
\end{figure}

Now we need a standard notion (see [6])
$\mathrm{ord}_ml$ that means the least integer number
$x$ s.t. $l^x=1 (\mathrm{mod}\ m)$.

\begin{lem}\label{fth.6}
$n$ odd $\implies (\forall
i\in\{0,\dots,n\!-\!1\})\exists!j\in\{0,\dots,n\!-\!1\})
[\natural(j)=i]$. So cycles of the mapping $\natural$
consist a partition of the set $\{0,\dots,n\!-\!1\}$.
One of the cycles is $\{0\}$. At last, the lengths of
any of these cycles divide the number
$\mathrm{ord}_n2$.
\end{lem}

{\bf Proof.} First let's define set
$M(s)=\{\natural^i(s)|i=0,1,2,\dots\}$ starting from
arbitrary position number $s\in\{0,\dots,n-1\}$ of
$X$. It's clear, all these sets are finite. For $s=0$
we have $\natural(0)=0$ and so $M(0)=\{0\}$. Set
$M(s),s>0,$ is a cycle iff there exists $i$ s.t.
$\natural^i(s)=s$. To prove that let's rewrite
elements of $M(s)$ in the natural order. We as before
will denote by $s$ the first element of this list.

One can note in first, that $i\neq2i(\mathrm{mod}\
n),i<n,n$ is odd. This mean that no other loop (cycle
with length 1) exists. Then, supposing the contrary
(i.e. $M(s)$ is not a cycle) we must conclude that for
some $m>0$ there exist at least two
$\natural$-prototypes $a,b\in\{0,\dots,n-1\},a\neq b,$
of $c=\natural^m(s)$, one of them, say $a$, is
$\natural^{m-1}(s)$. Because $2b=c(\mathrm{mod}\ n)$
and $b,c<n$ then $2b-c=n$. However because $c>s$ the
both of $2b,c$ are even numbers whereas $n$ is odd.
The contradiction proves that $M(s),0\notin M(s),$ is
a cycle because the unique opportunity is $m=0$, since
then $s$ can be odd. Moreover we conclude that the
least element of $M(s)$ in case $M(s)\neq\{0\}$ must
be odd number.

Now, because $2^{\mathrm{ord}_n2}=1(\mathrm{mod} n)$
the sequence $(2^is(\mathrm{mod} n))_{i=0,1,2,\dots}$
is is built from the recurrent pieces of length
$\mathrm{ord}_n2$. In other words, $|M(s)|$ divides
$\mathrm{ord}_n2.\ \ \ \Box$

How many of the cycles of kind $M(s)$ exist? At least
two in not trivial case $n>1$. The
example\ref{ex_ind_1} shows 5 cycles for $n=33$. In
this case $\mathrm{ord}_n2=10$.

On the other side when we pass from the "index
projection" $\natural$ to the baker transformation
$\frak b$ we must keep in mind that the series of
cycles in the set $\{\frak b^i*X|i=0,1,\dots\},|X|$ is
odd, could be different because of the distribution of
ones in $X$. Yet, it's obvious that the lengthes of
these cycles divide lengthes of the correspondent
cycles for $\natural$ and so must divide the same
number $\mathrm{ord}_{|X|}2$. However we do not need
the last statement to state the next corollary because
by definition the baker transformation,
proposition~\ref{prop_ab_bak_and in-bak}, and from the
previous lemma.

\begin{cor}\label{col.5}
If $n$ is odd then $\forall X\in\frak B^n[\frak
b^{\mathrm{ord}_n2}*X=X]$.
\end{cor}

The cases of even $n$ can be reduced to the previous
cases in the following way.

\begin{lem}\label{lem6}
Let $|X|=n=2n'$. Then $\frak
b*X=[y_0,0,y_1,\dots,y_{n'},0]$ for some boolean
numbers $y_0,\dots,y_{n'}$ and $\frak
b^2*X=[z_0,0,\dots,z_{n'},0]$ where
$[z_0,\dots,z_{n'}]=\frak b*[y_0,\dots,y_{n'}]$.
\end{lem}

{\bf Proof.} First of all, by the
definition~(\ref{baker_df}) only components of $\frak
b*X$ with numbers $2i(\mathrm{mod}\
n),i=\overline{0,n-1}$ can be non-zero. Taking in
account that all the numbers $2i(\mathrm{mod}\
n),i=\overline{0,n-1}$ are even (because $n$ is even),
we come to $\frak b*X=[y_0,0,y_1,\dots,y_{n'},0]$.

Then, when we apply $\frak b$ to
$[y_0,0,y_1,\dots,y_{n'},0]$ its odd components play
no role in forming components of the result, because
all they are equal to 0. It's important also that
indices $\natural(2i(\mathrm{mod}\ n))$ are even. So
the result $\frak b*[y_0,0,y_1,\dots,y_{n'},0]$ can be
write as $[z_0,0,\dots,z_{n'},0]$.

It remains only to note that we can shorten vectors
$[y_0,0,y_1,\dots,y_{n'},0]$, $[z_0,0,\dots,z_{n'},0]$
to $[y_0,y_1,\dots,y_{n'}]$, $[z_0,z_1,\dots,z_{n'}]$
because, as we said, the components of
$[y_0,0,y_1,\dots,y_{n'},0]$, $[z_0,0,\dots,z_{n'},0]$
with odd numbers play no role, and if $k=2k',t=2t'$
then $2k=i\ (\mathrm{mod}\ n)\ \&\ 2t=i\
(\mathrm{mod}\ n)\iff 2k'=\frac{i}{2}\ (\mathrm{mod}\
n')\ \&\ 2t'=\frac{i}{2}\ (\mathrm{mod}\ n').\ \ \
\Box$

Since it is possible to iterate the reduction of
lemma~\ref{lem6} we come to the next results.

Let $\iota(a,b)$ be the maximal degree of $a$ that
divides $b$ without a remainder and let's call any
boolean vector $X$ as {\sl $\frak b$-swept or
baker-swept} if $\forall j[0<j<|X|\ \&\
2^{\iota(2,|X|)}\nmid j\implies X_j=0]$. The prompt
consequence of the previous lemma is

\begin{cor}\label{col.7}
$\forall X[\frak b^{\iota(2,|X|)}*X$ is $\frak
b$-swept$]$.
\end{cor}

Now let us define a function $\frak c(n)$:
\begin{center}
$\frak c(n)=
\begin{cases}
\mathrm{ord}_{(n/2^{\iota(2,n)})}2,&
\text{ if } \frac{n}{2^{\iota(2,n)}}>1,\\
1,&  \text{ otherwise}.
\end{cases}$
\end{center}

\begin{theo}\label{th.7} $\forall n>0\forall
X\in\frak B^n[\frak b^{\iota(2,n)+\frak c(n)}*X= \frak
b^{\iota(2,n)}*X]$.
\end{theo}

{\bf Proof.} This results from lemma~\ref{lem6} about
reduction, previous corollary~\ref{col.6}, and
lemma~\ref{fth.6}. $\ \ \ \Box$

\begin{ex}\label{ex_ind_2}
{\rm The fig.~\ref{fig12} illustrate the theorem in
case $n=2^3\cdot33=264$. The diagram contains 5 cycles
of lengths the same as in case $n=33$. However, every
vertex of every cycle is the root of the same tree of
height 3.}
\end{ex}

\setlength{\unitlength}{.1cm}
\begin{figure}[htbp]
\begin{picture}(180,30)
\put(10,0){\includegraphics[width=10cm]{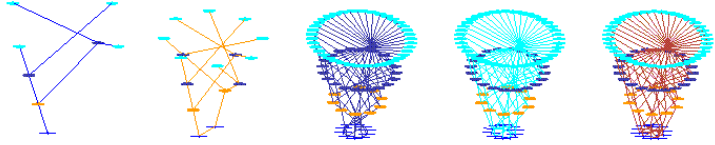}}
\end{picture}
\caption{The diagram of the mapping $\natural$ for
$n=264$.}\label{fig12}
\end{figure}

We call the number $\frak c(n)$ {\sl critical}. The
behavior of the function is shown by Fig.\ref{fig9}
for $n\le200$ (compare with fig. 4 from [6, p.85] and
note, please, that their picture is the graph of
$\mathrm{ord}_n2$ built only for odd numbers $n$).

\setlength{\unitlength}{.1cm}
\begin{figure}[htbp]

\begin{picture}(100,80)
\put(10,0){\includegraphics[width=8cm]{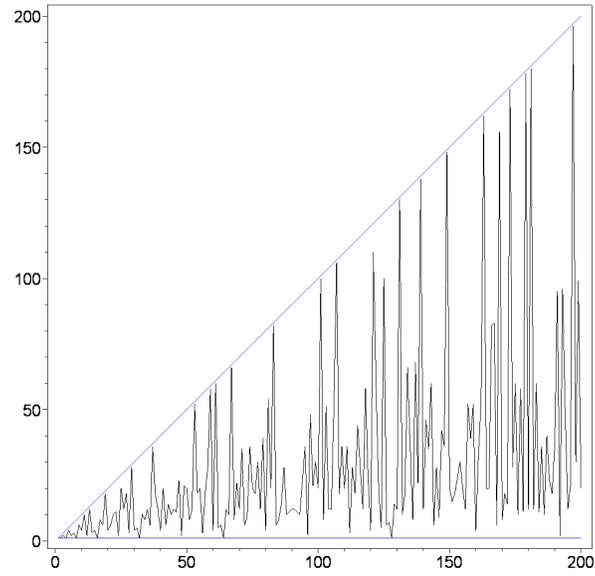}}
\end{picture}
\caption{Function $\frak c(n)$ on segment
$n=\overline{1,200}$.}\label{fig9}
\end{figure}

\begin{cor}~\label{entropy}
The numbers of the series
$$
\mathrm{rank}(C(X)),\mathrm{rank}(C(\frak b*X)),
\mathrm{rank}(C(\frak b^2*X)),\dots
$$
doesn't increase and their minimum equals to
$\mathrm{rank}(C(\frak b^{\iota(2,|X|)}*X))$.
\end{cor}

{\bf Proof.} Indeed, as we saw above $\frak
b^{\iota(2,|X|)}*X$ belongs to the period (cycle) of
the series $X,\frak b*X,\frak b^2*X,\dots,...$. Also
according to theorem~\ref{conserv_pr}(ii) the series
doesn't increase. $\ \ \ \Box$

The kneading ability of $\frak b$ is shown on the Fig.
\ref{fig10}. There $t$ is number of iteration of
$\frak b$. The beginning is subsegment $[.1,.2]$ of
the segment $[0,1]$ consisting of the numbers of kind
$x=\sum_{i=1}^{13}a_i2^{-i}$ where coefficients
$a_i,i=\overline{1,13}$ are boolean and $.1\le
x\le.2$.

\setlength{\unitlength}{.1cm}
\begin{figure}[htbp]
\begin{picture}(100,90)
\put(0,0){\includegraphics[width=9cm]{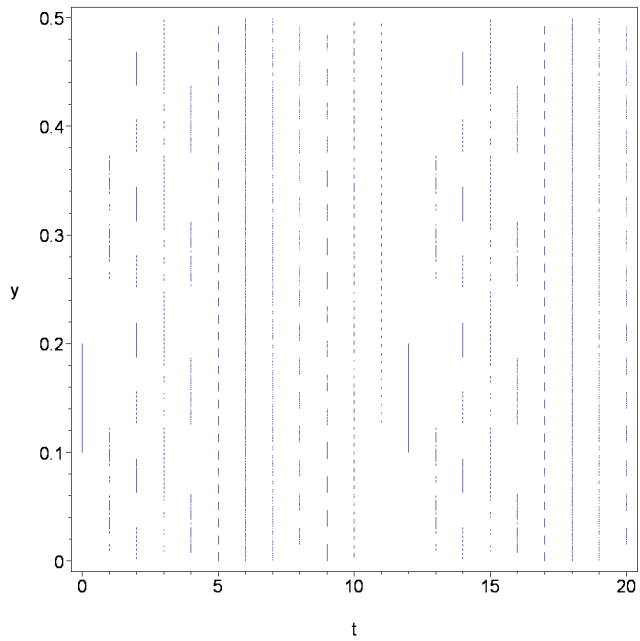}}
\end{picture}
\caption{Pseudo-chaotic behavior $\frak b$ on segment
$[0,1]$ for odd $n,n=13$.}\label{fig10}
\end{figure}

Despite of this overt kneading, the period of length
12 is here. This is because of the finiteness of the
mixed set: its power is $2^{13}$ or 8,192.

In case when set consists the rational numbers
representing boolean sequences of length $2^m$ the
picture is quite different as the fig.~\ref{fig10a}
shows.

\setlength{\unitlength}{.1cm}
\begin{figure}[htbp]
\begin{picture}(100,90)
\put(0,0){\includegraphics[width=9cm]{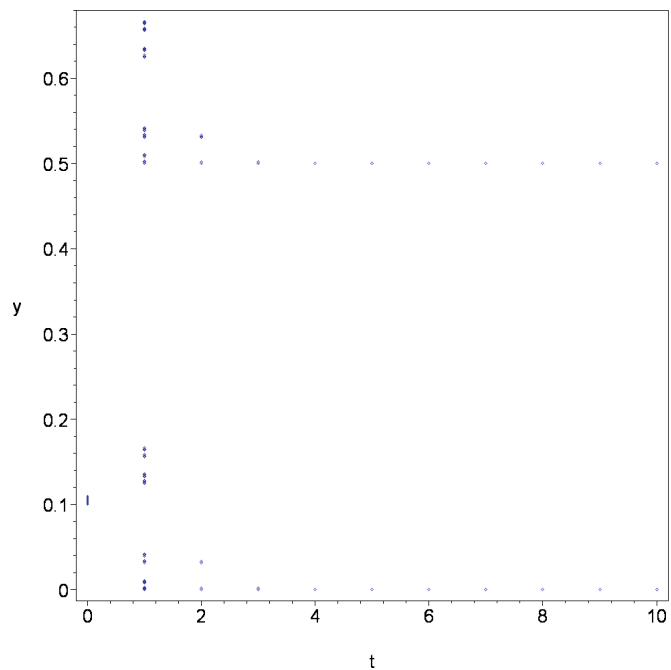}}
\end{picture}
\caption{Behavior $\frak b$ on segment $[0,1]$ for
$n=2^4$.}\label{fig10a}
\end{figure}
 \ \\

\section{Diagrams of $\frak b$ and its index projection
$\natural$}

As we saw the baker transformation is the useful tool
and therefore it's worth to gain more information
about it as well as about its index projection
$\natural$. The task become easier because we deal
actually with the same transformations (as functions)
for any space $\frak B^n$. On the other side for given
a number $n$ the diagrams of these mappings tell us
something about all rules of additive automata on the
cylinder of size $n$.

Clearly, $\natural$ acts on set of power $n$ whereas
space of $\frak b$ has power $2^n$. Therefore, for a
fixed $n$ the $\natural$-diagram presents
$rule$-independent data; whereas $\frak b$-diagram
differs individual rules $X\in\frak B^n$. The problem
arises, how to decode the information hidden in the
diagram. In this paragraph we, in particular, give
examples of the decoding.\\

\subsection{$\frak b$ as a linear operator
in $\frak B^n$}
 \ \\
\begin{lem}\label{lem7} For every $n$ the discrete baker
transformation acts as a linear operator on $\frak
B^n$.
\end{lem}

{\bf Proof.} It's enough to check:
\begin{enumerate}
\item $\frak b*0\cdot X=0\cdot\frak b*X$.

\item $\frak b*(X\oplus Y)=\frak b*X\oplus\frak b*Y$.
\end{enumerate}

Point (1) follows $\frak b*\vec0=\vec0$ that is the
prompt consequence of the definition~(\ref{baker_df}).

Point (2) also is true because of
\begin{eqnarray*}
(\frak b*(X\oplus Y))_i&=&\underset{j:\
2j=i(\mathrm{mod}\ n)}{\bigoplus}(X\oplus Y)_j=\\
&=&\underset{j:\ 2j=i(\mathrm{mod}\
n)}{\bigoplus}(X)_j\ \oplus\ \ \underset{j:\
2j=i(\mathrm{mod}\ n)}{\bigoplus}(Y)_j=\\
&=&(\frak b*X)_i\oplus(\frak b*Y)_i.
\end{eqnarray*}
$\ \ \ \Box$

Because $\frak b$ is the linear operator it can be
represented by matrix $\mathcal B_n$ for given a
dimension $n$. As for operators $\hat X$ we accept
that operator and the corresponding matrix act on
vectors (strings) from left side. The matrixes for
even $n$ essentially differs from the case of odd $n$.

%%%%matrixes B_8,B_9 are excluded

In general for $i,j=\overline{0,n-1}$
\begin{center}
$(\mathcal B_n)_{i,j}=
\begin{cases}
1,& \text{ if } i=\mathrm{rem}(2j,n),\\
0,&  \text{ else}.
\end{cases}$
\end{center}
(Here, as before, $\mathrm{rem}(a,b)$ denotes
remainder of division $a$ with $b$.)

So despite of additivity $\frak b$ differs from
cellular automata as operators.

Now we present examples of $\frak b$-diagrams.

\begin{figure}[htbp]
\begin{center}
\includegraphics[width=4cm]{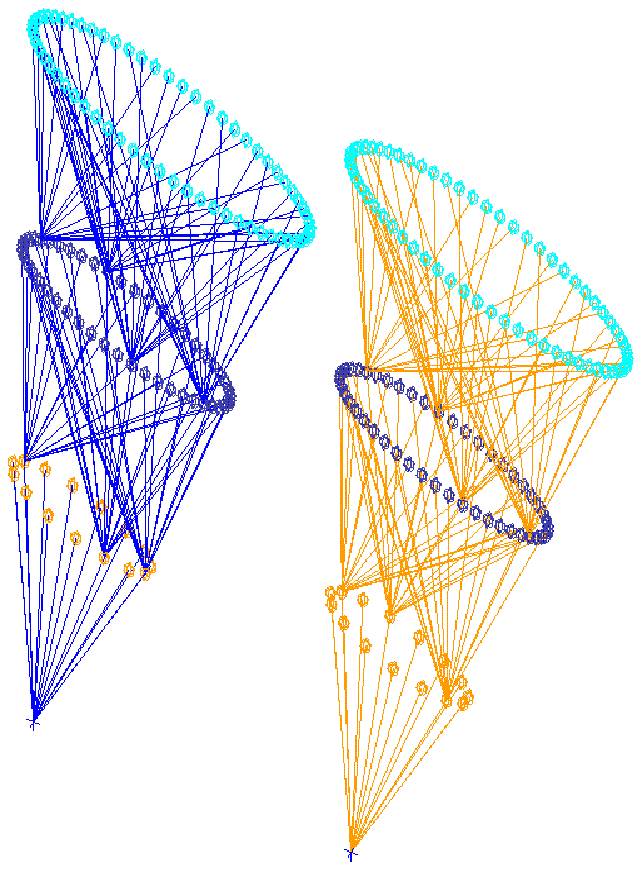}
\end{center}
\caption{The diagram of $\frak b$ in $\frak
B^8$.}\label{fig7}
\end{figure}

\begin{ex}\label{ex.7}
{\rm Generally speaking, when $n=2^k$ the behavior of
$\frak b$ in $\frak B^n$ must be sufficiently
predictable in view of theorem~\ref{th.7}. We will
discuss this later. Fig.~\ref{fig7} represents the
diagram of behavior $\frak b$ in $\frak B^8$.}
\end{ex}

%\newpage

\begin{ex}\label{ex.8}
{\rm The diagram for $\frak B^9$ is the collections of
8 cycles with length 1, 12 cycles with length 2, 8
cycles with length 3, and 76 cycles with length 6.}
\end{ex}

\begin{figure}[htbp]
\begin{picture}(160,60)
\put(10,0){\includegraphics[width=4cm]{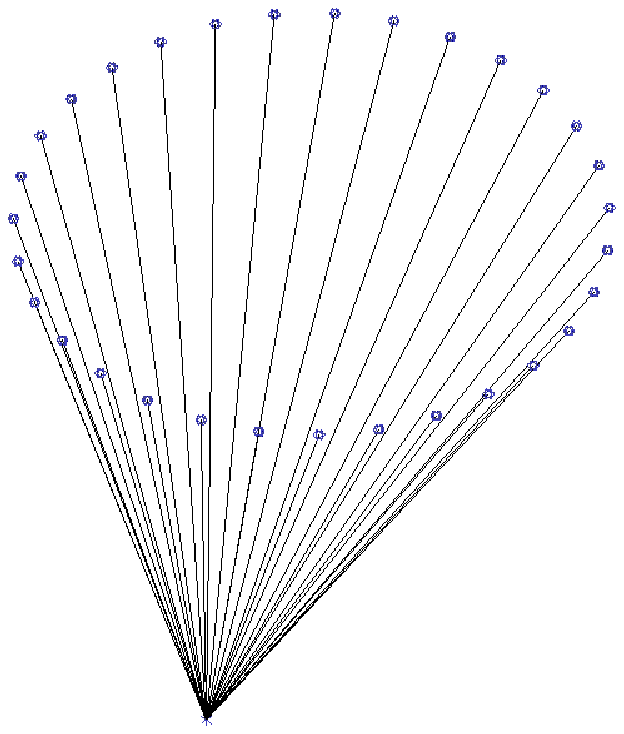}}
\put(90,0){\includegraphics[width=4cm]{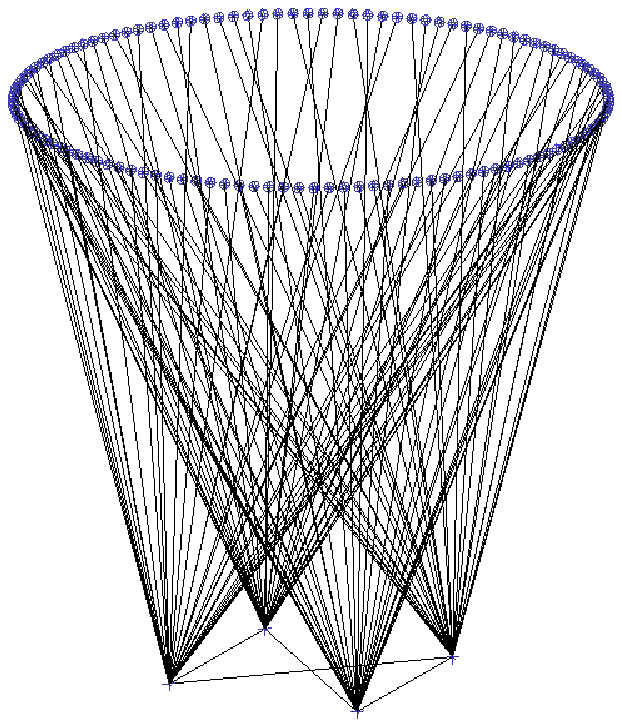}}
\put(50,0){\includegraphics[width=4cm]{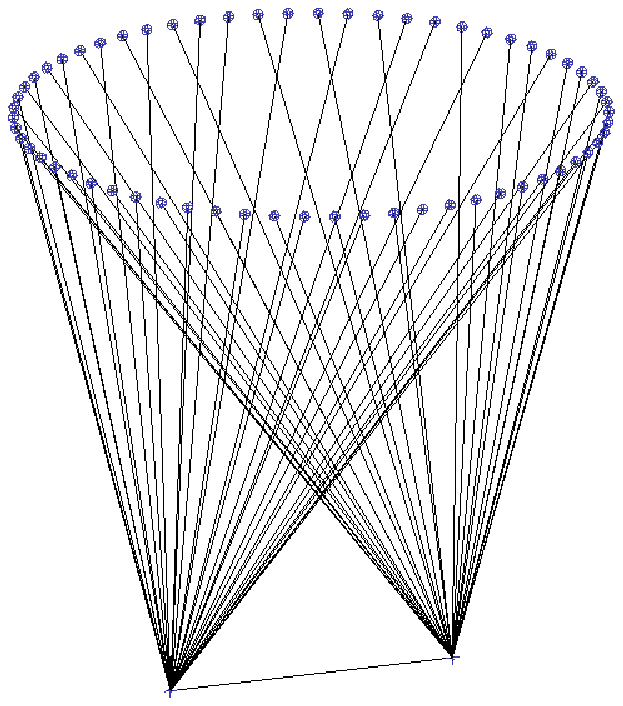}}
\end{picture}
\caption{Three kinds of connectivity components of the
diagram of $\frak b$ in $\frak B^{10}$.}\label{fig8}
\end{figure}

\begin{ex}\label{ex.9}
{\rm Fig. \ref{fig8} presents three kinds of
connectivity components of the baker diagram in $\frak
B^{10}$. The complete diagram includes 4 basins of
cycles of length 1, 2 basins of cycles of length 2 and
6 basins with cycles of length 4.}
\end{ex}

\ \\

\subsection{$\natural$ as a homomorphism of
finite cyclic groups}
 \ \\
We denote by $\mathbb C_n$ the additive cyclic group
$\langle \{0,1,\dots,n-1\},+_n\rangle$ where $+_n$ is
the addition w.r.t. modulo $n$ on $\{0,\dots,n-1\}$.

\begin{lem}\label{lem8}
$\natural$ is a homomorphism of $\mathbb C_n$ on some
its subgroup which is an isomorphism if $n$ is odd.
\end{lem}

{\bf Proof.} By the definition~(\ref{ind_baker_df})
$\natural(i)=\mathrm{rem}(2i,n)$; so $\natural(0)=0$
and $\natural(i+_nj)=\mathrm{rem}(2i+_n2j,n)=
\mathrm{rem}(2i,n)+_n\mathrm{rem}(2j,n)=
\natural(i)+_n \natural(j). \ \ \Box$

Thus the proposition~\ref{prop_ab_kernels} is
applicable here too, and this explains the fact that
any diagram of $\natural$ consists several cycles, one
of them is $\{0\}$; and every vertex of every cycle is
the root of the same (within isomorphism, of course)
tree. Examples\ref{ex_ind_1},\ref{ex_ind_2}
demonstrate the diagrams of $\natural$ in $\mathbb
C_{33},\mathbb C_{264}$ correspondingly. The big
advantage of $\natural$ is that it acts in small space
$\mathbb C_n$ comparing with $\frak B^n$ for $\frak b$
or with $2^n$ STD of $\frak A(n,X)$ when $X$ run
$\frak B^n$. Indeed, it's impossible even to compute
(to not say about visualizing) diagram for $\frak B$
in $\frak B^{264}$.

The other nice feature of $\natural$-diagram is that
shrinking factor is 2 for every level of the diagram
apart from the level of cyclic vertexes.

\section{Reading baker diagrams}

Here we start from upper estimates for height and
cycle lengths. Then we show that these estimates are
not improvable. And then we apply the results for
decoding information given by baker and index-baker
diagrams to characterize in whole the system of all
ACA of the cylinder of size $n$.

Firstly, due to the conservation principle for $\frak
b$ all rules of the same connectivity component of the
baker diagram in $\frak B^n$ have the same determinant
value. This, for example, means that if one rule of a
connectivity component has $h^{\star}=0$ then the same
is true for all others.

We say that $X\in\frak B^n$ {\sl belong to the cycle}
of the baker transformation if there exists a number
$t$ s.t. $\frak b^t*X=X$.

\begin{lem}\label{lem.5}
 If $X$ belongs to a cycle of the
baker transformation diagram in $\frak B^n$ then
$h^{\star}(n,X)\le1$.
\end{lem}

{\bf Proof.} By the condition and theorem~\ref{th.5}
we have $\hat X^{2^{it}}=\hat X$ for fixed a number
$t$ and any $i,i>0$. Therefore $\ker(\hat X)=\ker(\hat
X^2)$ (see proposition~\ref{prop_ab_kernels}(5)). This
means that images of the operators $\hat X,\hat X^2$
also are the same and $\ker(\hat
X)\cap\mathrm{Im}(\hat X)=
\{\mathbf0\}$. Hence two opportunity remain:\\
(i) $\ker(\hat X)=\{\mathbf 0\}$. If so then obviously
$h^{\star}(n,X)=0$.\\
(ii) $\ker(\hat X)\neq\{\mathbf 0\}$.
Then\footnote{$\hat X^0$ is the identical operator.}
$\mathrm{Im}(\hat X^0)\supsetneq\mathrm{Im}(\hat X)$
and the correspondent shrinking factor $|\ker(\hat
X)|>1$. So all states from $\mathrm{Im}(\hat X)$ are
in cycles and states from $\mathrm{Im}(\hat
X^0)\setminus\mathrm{Im}(\hat X)$ are not. So,
$h^{\star}=1$.
 $\ \ \ \Box$

\begin{cor}\label{col.3} If $n$ is odd then
$h^{\star}(n,X)=1-\det_2(C(X))$.\footnote{As we agreed
in the beginning $\det_2(A)$ is the determinant modulo
2 of a matrix $A$.}
\end{cor}

{\bf Proof.} This follows from the previous lemma and
theorem~\ref{th.7} above. $\ \ \ \Box$

\begin{cor}\label{col.4} (Lemma 12[1].)
If $\kappa(X)$ is even and $n$ is odd then
$h^{\star}(n,X)=1$.
\end{cor}

{\bf Proof.} Indeed, if the parity of $X$ is 0 then
$\det_2(C(X))=0$ because sum of all column modulo 2
equals 0 and therefore the rank (in the boolean field)
of matrix $C(X)$ isn't $|X|$. $\ \ \Box$

\ \\

\subsection{Upper height estimations}
 \ \\
Let's start from height estimation. To estimate
$h^{\star}$ let's note that the least effective value
of shrinking factor can be 2 for all levels apart from
level 0 consisting cyclic states only and where no
shrinking exist. Indeed, if there exist non-cyclic
states then the kernel of the rule has power bigger 1
and in-degree of any state is 0 for the dangling
states and bigger 1 for the other, i.e. for states
that have prototypes (see the
proposition~\ref{prop_ab_kernels}).

Therefore we can write
$t+t(\sum_{i=1}^{h^{\star}}2^{i-1})=|\frak B^n|$,
where $t$ is the number of cyclic states. The right
part equals to $t2^{h^{\star}}$. From here, setting
the minimal value for $t,t=1,$ we get
$h^{\star}(n,X)\le\log_2|\frak B^n|=n$. This estimate
obtained by "bare hands" can be improved as the
following.

\begin{theo}\label{th8}
$\forall n(\forall X\in\frak
B^n)[h^{\star}(n,X)\le2^{\iota(2,n)}]$.
\end{theo}

{\bf Proof.} It is an easy consequence of the
theorems~\ref{th.7} and \ref{th.5}. Indeed, from the
theorem~\ref{th.7} it follows that $\frak
b^{\iota(2,n)}*X$ is a cyclic state. And then from
theorem~\ref{th.5} we get that the superposition $\hat
X^{2^{\iota(2,n)}}$ transfers any state $s$ into a
cycle. $\ \ \ \Box$

The number $\iota(2,n)$ is easy computable on $n$ and
in the same time is the maximal height of index-baker
(and baker diagrams too). On the other side, the baker
diagram distributes rules $X\in\frak B^n$ on levels
with the same distance from attractors and therefore
contains additional information about rules in
comparison with index-baker diagram. So we can improve
the estimate above for given a rule $X$. Let denote by
$H(n,X)$ the distance from $X$ to the closest cyclic
vertex in the baker diagram of $\frak B^n$.\footnote{
Generally speaking to compute $H(n,X)$ one don't need
to compute complete baker diagram for $n$. However
sometimes the complexity can be almost the same.}

\begin{theo}\label{th.9}
$(\forall X\in\frak
B^n)[h^{\star}(n,X)\le2^{H(n,X)}]$.
\end{theo}

{\bf Proof.} The proof actually the same as for the
previous theorem because we use the analogous fact: if
$H(n,X)=m$ then $\frak b^m*X$ belong to a cycle. That
mean $\hat X^{2^m}*s$ (see theorem~\ref{th.5}) is
cyclic state for any initial state $s\in\frak B^n$.
The last mean $h^{\star}(n,X)\le2^m.\ \ \Box$

For example, in case $n=10$ this estimation looks as
$h^{\star}\le2$. And the value 2 is reached for
$X=[1,1,1,1,1]$ as well as for $X'=[1,1,0,1,0,1]$, see
Fig. \ref{fig.11}.

\begin{ex}\label{ex10}
{\rm The next Fig. \ref{fig.11} shows STD for
$n=10,X=[1,1,0,1,0,1]$. Here $h^{\star}=2$ as it is
estimated by theorem~\ref{th.9}. So the upper estimate
was reached for this rule $X$. The "bare hands"
estimate 10 is essentially bigger.}
\end{ex}

\begin{figure}[htbp]
\begin{center}
\includegraphics[width=8cm]{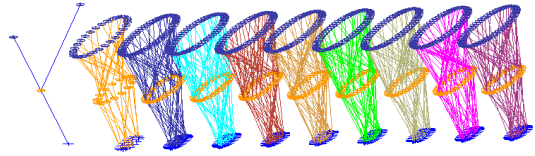}
\end{center}
\caption{STD for
$n=10,X=[1,1,0,1,0,1]$.}\label{fig.11}
\end{figure}

\subsection{Cycle length estimation}
\ \\
\begin{theo}\label{th10}
For any given $n$ and any $X\in\frak B^n$ lengths of
cycles in STD of $\frak A(n,X)$ must divide number
$c^{\star}(n)=2^{\iota(2,n)}(2^{\frak c(n)} -1)$.
\end{theo}

{\bf Proof.} According to
theorems~\ref{th.7},\ref{th.5} we have $\forall s[\hat
X^{\exp_2(\iota(2,n)+\frak c(n))}*s=\hat
X^{\exp_2(\iota(2,n))}*s]$ where $\exp_ab=a^b$.
Because state $s'=\hat X^{\exp_2(\iota(2,n))}*s$
belongs to a cycle we can inverse $\hat
X^{\exp_2(\iota(2,n))}$ on $s'$ getting another state
$s''$ from the cycle $s''=\hat
X^{-\exp_2(\iota(2,n))}*s'$, see fig.~\ref{pic_cyc}.
Hence the relation is true
\begin{equation}\label{s-s''}
\hat X^{\exp_2(\iota(2,n)(\frak c(n)-1)}*s=s''.
\end{equation}
Then we can replace $s$ with $s''$ in~(\ref{s-s''})
and having this done we come to the relation for the
element $c''$ of the cycle
\begin{equation*}
\hat X^{\exp_2(\iota(2,n))(\exp_2(\frak
c(n))-1)}*s''=s''.
\end{equation*}
This just mean that the number $c^{\star}$ must be
divided without any remainder by the length of the
cycle. Also since $s$ was any initial state, it is
possible to say this about any cycle of the diagram
for the rule $X$.
 $\ \ \ \Box$

\setlength{\unitlength}{.1cm}
\begin{figure}
\begin{picture}(80,60)

%\begin{picture}(80,60)
\put(10,10){\includegraphics[width=5cm]{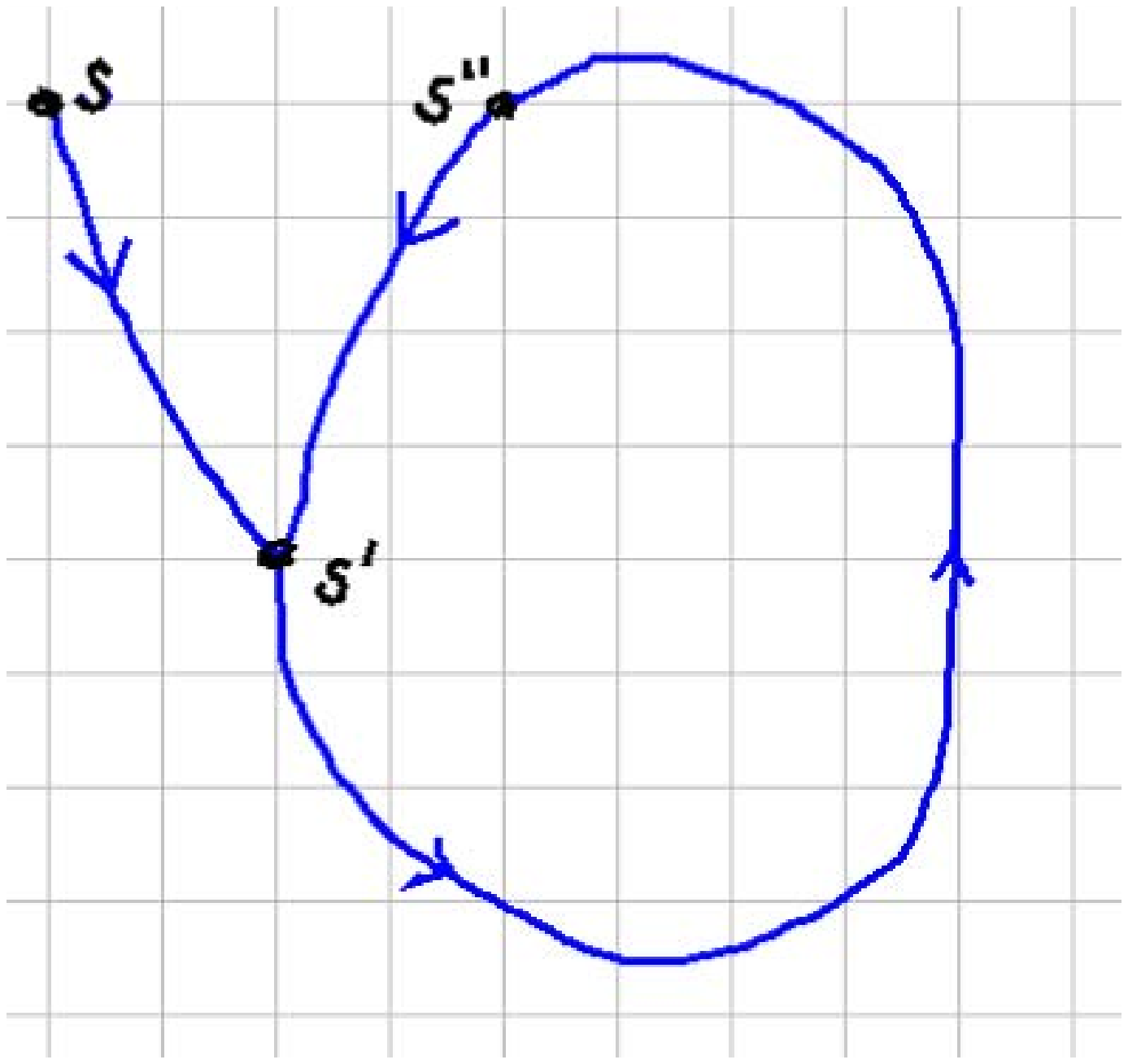}}
%\end{picture}

 \put(0,40){$\hat X^{2^{\iota(2,n)}}$}
  \put(30,40){$\hat X^{-2^{\iota(2,n)}}$}
 \put(12,45){\vector(2,-3){8}}
 \put(25,35){\vector(1,2){6}}
% \put(-82,40){$\hat X^{2^{\iota(2,n)}}$}
%  \put(-50,40){$\hat X^{-2^{\iota(2,n)}}$}
%\put(-40,32){\circle{20}}
 %\put(-41,24){$\succ$}
% \put(-70,45){\vector(2,-3){8}}
% \put(-55,35){\vector(1,2){6}}
\end{picture}

\vspace{-1.0cm} \caption{To the proof of the
theorem~\ref{th10}.}\label{pic_cyc}
\end{figure}

In the diagram of $\frak A(8,[1,1,1])$ there are 2
cycles with length 1, 1 cycle of length 2, 3 - of
length 4, and 30 cycles of length 8. So all possible
for $n=8$ lengths are realized because
$c^{\star}(8)=2^{3}(2^{1} -1)=8$.

However this happens not every time. Dimension $n$
constrains scope of possible lengths of cycles. The
real collection of cycle lengths in STD of $\frak
A(n,X)$ for a fixed $n$ varies with $X$.

If $n=10$ then $\frak c(10)=4,2^4-1=15$, so lengths
can run set of dividers of 30. For STD of $\frak
A(10,[1,1,0,1,0,1])$ (see example before, Fig.
\ref{fig.11}) we have one fixed point $\mathbf 0$, one
cycle of length 15 and 8 cycles of length 30;
$h^{\star}=2$.

However for the rule $[1,1,1,1]$ with $h^{\star}=1$
there are 24 cycles of length 10, three of length 5,
and one fixed point; whereas for rule
$[0,0,0,0,0,0,1,1,1,1]$ with $h^{\star}=1$ there exist
16 cycles of length 1 and 120 cycles of length 2.

At last for rule $[1, 0, 0, 1, 0, 1, 1, 0, 1, 1]$ we
have $h^{\star}=2$ and 1 cycle of length 1, 5 cycles
of length 3, and 40 cycles of length 6.

One more example: $\frak c(9)=6,\ \iota(2,9)=0$, so
any cycle of $\frak A(9,X)$ for any $X$ must be a
divider of $2^6-1=63$. STD for $X=[1,1,0,1,0,0,1,0,1]$
has cycles of lengths 7, STD for $X=[1,1,0,1,0,1,1]$
has 6 cycles of length 21, and STD for $[1,1]$
contains 4 cycles of length 63, but STD for
$X=[1,0,0,1,0,1,1]$ - has no cycles with lengths that
are multiple of 7.

As far we use rule independent computationally easy
information about cycle lengths that we get from
actually index projection $\natural$ of the baker
transformation. Again the baker diagram provides, in
general, more accurate estimates depending on concrete
rules $X$. True, this information is more costly in
the sense of computability.

As before we use denotation $H(n,X)$ for the height of
the rule $X$ in the diagram of $\frak b$ in $\frak
B^n$. In addition, let's denote by $\frak C(n,X)$ the
length of cycle in basin of which $X$ is for the
diagram.

\begin{theo}\label{th11}
For given $n$ and $X\in\frak B^n$ lengths of cycles in
STD of $\frak A(n,X)$ must divide number
$C^{\star}(n,X)=2^{H(n,X)}(2^{\frak C(n,X)} -1)$.
\end{theo}

The proof of the theorem repeats the proof in general
case with the natural replacement $\iota(2,n)$ and
$\frak c(n)$ with the numbers $H(n,X)$ and $\frak
C(n,X)$ that are specified for the considered rule
$X$. (It is worth to note that $2^{H(n,X)}(2^{\frak
C(n,X)} -1)|2^{\iota(2,n)}(2^{\frak c(n)} -1)$, so no
contradiction exists between the last two theorems.)

Now let's apply the last result to the case $n=10$.

The rules $X=[1,1,0,1,0,1]$ and $X=[1,1,1,1]$ belong
to basins of cycles with length 4 in the baker
diagram. The both of them have $H=1$. So the last
theorem tell us nothing of new.

However the rule $X=[0,0,0,0,0,0,1,1,1,1]$ with
$H(10,X)=1$ belongs to basin of cycle with length 1.
Therefore only cycles of lengths 1 and 2 can occur in
STD of $\frak A(10,X)$. And indeed, our computation
brings 16 cycles of length 1, 120 cycles of length 2.

Then, rule $[1, 0, 0, 1, 0, 1, 1, 0, 1, 1]$ occurs in
a basin of the cycle of length 2 in the baker diagram.
Therefore in accordance with theorem~\ref{th11} we
have 1 cycle of length 1, 5 cycles of length 3, and 40
cycles of length 6.

At last, rule $X=[1, 0, 0, 0, 1, 0, 1, 0, 1, 0]$
belongs to the cycle of length 4 in the baker diagram.
So $H(10,X)=0,\frak C(10,X)=4$. The theorem~\ref{th11}
states that the cycles of $\frak A(10,X)$ must divide
the number $2^0(2^4-1)$=15. And indeed, according to
computations STD for $\frak A(10,X)$ contains 1 fixed
point and 51 cycles of length 5. Note,
$h^{\star}(10,X)=1$.

The data are represented by the table, see
fig.~\ref{pic_cyc}. We recall, that $H^{\star}=1,\
c^{\star}=30$; so the supposed on the base of index
baker diagram collection of cycle lengths is
1,2,3,5,6,10,15,30. Knowledge of the position of rule
$X$ (i.e. $H(10,X)$ and the power $\frak C(10,X)$ of
the attractor (=cycle) whose basin contains $X$) in
the baker diagram often provide more exact estimation.

\footnotesize
\begin{figure}[htbp]
\begin{tabular}{|c||c|c|c|c|c|c|}
 \hline
 N&$X$&attractor&$X$-height&$C^{\star}(10,X)$
 &cycle&$h^{\star}(10,X)$\\
 && power $\frak C(10,X)$&$H(10,X)$&& lengths&\\
 \hline\hline
   1 &$[0000001111]$&1&1&2&1,2&1\\
    \hline
      2 &$[1001011011]$&2&1&6&1,3,6&2\\
    \hline
      3 &$[1000101010]$&4&0&15&1,5&1\\
    \hline
    4&$[1101010000]$&4&1&30&1,15,30&2\\
  \hline
  5&$[1111000000]$&4&1&30&1,5,10&1\\
  \hline
\end{tabular}

\caption{Cycle lengths for different rules,
$n=10$.}\label{pic_cyc}
\end{figure}
\normalsize

As another example of the theorem application we sum
up the result relating to odd number $n=9$ in the
table of fig.~\ref{pic_cyc1}. This time $H^{\star}=0,\
{\mathrm ord}_92=6,\ c^{\star}=2^6-1=63$; so the
supposed on the base of index baker diagram collection
of cycle lengths is 1,3,7,9,21,63. Now $H(9,X)=0$ but
the power $\frak C(10,X)$ of the attractor (=cycle)
whose basin contains $X$ in the baker diagram indeed
help getting more exact estimation, see the first 4
rows of the table.

\footnotesize
\begin{figure}[htbp]
\begin{tabular}{|c||c|c|c|c|c|}
 \hline
 N&$X$&attractor&$C^{\star}(9,X)$
 &cycle&$h^{\star}(9,X)$\\
 && power $\frak C(9,X)$&& lengths&\\
 \hline\hline
 1&$[000100100]$&1&1&1&1\\
  \hline
  2&$[101101101]$&2&3&1,3&1\\
   \hline
   3 &$[111000011]$&3&7&1,7&0\\
    \hline
      4 &$[110100101]$&3&7&1,7&1\\
    \hline
      5 &$[111000010]$&6&63&1,3,63&1\\
    \hline
          6 &$[110101100]$&6&63&1,21&1\\
\hline
      7 &$[100101100]$&6&63&1,3,9&1\\
\hline
\end{tabular}

\caption{Cycle lengths for different rules,
$n=9$.}\label{pic_cyc1}
\end{figure}
\normalsize

The comparison of the theoretical estimates and
experimental data done in these tables show that a
rule position in baker diagram doesn't provide
complete information about the spectrum of cycle
lengths and $h^{\star}$ value of STD for the rule.

This statement is confirmed by the following. Results
about lower estimations of $h^{\star}$ and maximal
cycle length would be of great interest. Yet, the
proposition:\\
\begin{equation}
2^{H(n,X)-1}(1-{\det}_2(n,X))\le
h^{\star}(n,X)\le2^{H(n,X)}(1-{\det}_2(n,X))
\end{equation}
looks pretty naturally and has many supporting it
examples. Nevertheless the example~\ref{ex11}
disproves it as the general statement.

\begin{ex}\label{ex11}
{\rm Let $n=12$. The rules
$X=[1,1,0,0,0,0,0,0,1,1,0,0]$ and
$Y=[1,0,1,0,1,0,1,0,0,0,0,0]$ are in the basin of the
rule $Z=[0,0,0,0,1,0,0,0,1,0,0,0]$. Moreover we have
$\frak b*X=Y,\frak b*Y=Z,frak b*Z=Z$. Despite
$H(12,X)=2>1=H(12,Y)>0=H(12,Z)$ for all these rules
$h^{\star}=1$. The difference takes place for cycle
spectra. $X$ comes with cycle lengths 1(4)\footnote{In
the round bracket we place the quantity of the cycles
with this length}, 2(6), and 4(60). For $Y$ there are
numbers 1(16), 2(120). And for $Z$ we have 1(). It is
a striking thing that $\ker\hat X=\ker\hat Y=\ker\hat
Z$. }
\end{ex}

So, not trivial lower estimations of $h^{\star}$
probably need more informative characteristics of a
rule than $H$. One of the appropriate tools is
obtaining below estimations by means of the special
imbedding one diagram into another the generating
rules of which are connected with each other by baker
transformation. We will not develop this idea here.
Instead further we check whether the language of DBT
equalities and inequalities is able to give us more.

The results about upper estimations tell us that when
we are moving from dangled vertexes of a baker diagram
to its attractor then upper bounds of $h^{\star}$ and
maximal lengths of cycles monotonically decrease. But
this is also true for the real maximal cycle lengths
and $h^{\star}$ ({\sl\bf
monotonicity principle}).\\

\section{Equalities and inequalities with $\frak b$
and $\boxtimes$}

This section continues the study how to interpret
baker diagrams; now we use equations with the
introduced operators.\\

\subsection{Equalities and inequalities with $\frak b$}
\ \\
Application of $\frak b$ to a rule $X$ in terms of a
baker diagram $G$ mean a passage from vertex $X$ to
the end of edge $(X,\frak b*X)$ of the diagram.
Therefore it's possible to express some relations and
substructures of the graph $G$. In this way some
conditions for rules can be set and solved. In
general, expressive power of the first-order language
with functions $\frak b,\natural$ on finite strings is
not a simple problem.

If we restrict ourselves with non-quantified formulas,
we come to systems of equalities and their negations.
The use for our theme of equalities and inequalities
with $\frak b$ can be seen clearly on the next
example.

\begin{lem}\label{lem9}
Three statement are equivalent:\\
(1) $\frak b^{\frak c(|X|)}*X=X$;\\
(2) ($X$ belongs to the cycle of the baker diagram in
$\frak B^{|X|}$);\\
(3) $\forall j[0<j<|X|\ \&\ 2^{\iota(2,|X|)}\nmid
j\implies X_j=0]$.
\end{lem}

{\bf Proof.} Firstly let $|X|$ be odd. Then, according
to theorem~\ref{th.7} (1) is true for all $\frak
B^{|X|}$ because $\iota(2,|X|)=0$; the baker diagram
consists of cycles only; and (3) is true trivially.

Now, let $|X|$ is even.

(1)$\implies$(2) because in general $(Y,\frak b*Y)$ is
the edge of the baker diagram and form (1) it follows
that starting from $X$ by means of $\frak c(|X|)$-edge
path we come back to $X$.

To show $(2)\implies(3)$ let on the contrary there be
a number $j,0<j<|X|,2^{\iota(2,|X|)}\nmid j,$ s.t.
$j$-th component of $X$ is not 0. Since $\frak
c(|X|)>0$ it's possible to pass to the result of any
finite applications of $\frak b$ to $X$ being inside
the cycle. By the definition of $\frak b$ $i$-th
application of $\frak b$ replaces all components of
the argument, components whose numbers are of kind
$2^{i-1}m,m$ odd, with 0 ("sweeps out" the
components). So, after $\iota(2,|X|)$ subsequent
applications of $\frak b$ to $X$ the condition (3)
must be true and never more any of these component can
be 1. If so, then $X$ can not belong to the cycle,
contradiction.

Now, let (3) is true. Then application $\frak b$ to
$X$ actually can be reduced to application of $\frak
b$ to the rule $X'\in\frak B^{|X|/\iota(2,|X|)}$ by
subsequent applying lemma\ref{lem6}. Since
$|X|/\iota(2,|X|)$ is odd number, it's clear that $X$
is in a cycle. It remains only to remind that length
of any cycle of the baker diagram divides the number
$\frak c(|X|)$.\footnote{One have not to mix cycles of
baker diagrams with cycles of automata STD. In the
last case see theorems~\ref{th10},\ref{th11}.} $\ \
\Box$

The condition (1) of the lemma is quite computable
within complexity $\mathcal O(|X|^3)$ since $\frak
c(|X|)\le|X|-1$ and the complexity of computing $\frak
b*X$ does not exceed $|X|^2$. However the condition
(3) can be checked for linear time relatively $|X|$.

Another remark is that despite the property of a rule
$X$ to be in a baker-cycle generally looks as $\exists
i[\frak b^i*X=X]$, in reality, as the lemma tell us,
this quantifier is bounded.

Recall, (corollary~\ref{col.7}), that we named a rule
$X$ {\sl $\frak b$-swept or baker-swept} if $\forall
j[0<j<|X|\ \&\ 2^{\iota(2,|X|)}\nmid j\implies
X_j=0]$. So the previous theorem states that $X$
belongs to a cycle of the baker diagram in $\frak
B^{|X|}$ iff $X$-is $\frak b$-compressed. Evidently,
the property to be baker-compressed is easily checked.
Therefore the next corollary has a sense.

\begin{cor}\label{col.6}
($\frak b^i*X$ is a $\frak b$-swept) $\implies
h^{\star}(|X|,X)\le2^i$.
\end{cor}

The proof is evident due the theorem~\ref{th.9} and
the previous lemma.

The sense of fixed point consideration for the baker
transformation becomes clear in view of the next
results.

Theorem~\ref{th.7} states that any rule $X$ satisfies
the equation of the kind $\frak b^q*X=\frak b^r*X$.
The exponents suggested in the theorem are common for
all the rules of length $n$. However particular rules
can satisfy also other equations of that kind. The
lower number $q,r$ the stronger restriction put on
$X$.

\begin{lem}~\label{lem10}
If for a rule $X$ there exists $q,r\in\mathbb N$ s.t.
$q>r>0\ \&\ \frak b^q*X=\frak b^r*X$ then
$h^{\star}(|X|,X)\le2^r$ and the least common multiple
of cycle lengths of STD for $\frak A(|X|,X)$ divides
the number $2^r(2^{q-r}-1)$.
\end{lem}

{\bf Proof.} The demonstration is in essence the same
as proof of theorems~\ref{th.9} and~\ref{th11}. $\ \
\Box$

The next result is a curios consequence of the lemma
because states a relation between fixed points $X$ of
$\frak b$ and fixed points of $\hat X$.

Here and below we denote the identical operator as
$\hat I$ and, naturally, $I=[1,0,\dots,0]$.

\begin{theo}~\label{the13} A rule $X$ is a fixed point
of $\frak b$ (i.e. $\frak b*X=X$) $\iff$ every
attractor of STD for $\frak A(|X|,X)$ is a fixed point
of $\hat X$ and ($\hat X=\hat I\ \vee \
h^{\star}(|X|,X)=1$).
\end{theo}

{\bf Proof.} $\Rightarrow$. As it follows from the
previous lemma for $q=1,r=0$ the all cycle lengths of
$\frak A(|X|,X)$ must divide $2-1=1$, i.e. the
attractor set of STD for $X$ consists of fixed points
only. Also $h^{\star}\le 1$.

Now, if $\det_2(|X|,X)=1$ then, since every state is a
fixed point of $\hat X$, we get $\hat X=\hat I$, i.e.
the identical operator. Otherwise, the determinant
equals to 0 and $h^{\star}=1$.

$\Leftarrow$. $\hat X=I$ means $X=[1,0,\dots,0]$.
Therefore $\frak b*X=X$. Now, let $h^{\star}(|X|,X)=1$
and every attractor of $\frak A(|X|,X)$ is a fixed
point. Then for any state $s$ the state $\hat X*s$ is
a fixed point. So $\forall s[\hat X^2*s=\hat X*s]$.
This means $\hat X^2=\hat X$ and therefore $\frak
b*X=X$.\footnote{It's easy to find $s$ s.t. $\hat
X*s\neq\hat Y*s$ if strings $X$ and $Y$ are not
equal.} $\ \ \ \Box$

When every attractor consists  a single state, the
quantity of the basins is a degree of 2 since it
equals to $\frac{2^{|X|}}{|\ker X|}$.

Solution of the equations of the considered type
doesn't present big difficulties at least in case when
the dimension of $X$ is given. For example, let $n=9$
and $X=[a_0,a_1,\dots,a_8]$. The equation looks as
$[a_0,a_1,a_2,a_3,a_4,a_5,a_6,a_7,a_8]=
[a_0,a_5,a_1,a_6,a_2,a_7,a_3,a_8,a_4]$. So it's not
problem to write the general solution:
$a_0=a,a_3=a_6=c,a_1=a_2=a_4=a_5=a_7=a_8=b$ or the
solution is $X_s=[a,b,b,c,b,b,c,b,b]$ where $a,b,c$
are arbitrary boolean numbers.

\begin{figure}[htbp]
\begin{center}
\includegraphics[width=6cm]{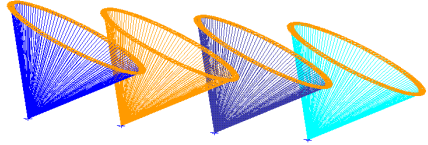}
\end{center}
\caption{The diagram for
$[0,1,1,0,1,1,0,1,1]$.}\label{fig.13}
\end{figure}

Let $a=0,b=1,c=0$. The diagram is represented on Fig.
\ref{fig.13}.

Another example of lemma~\ref{lem10}:

\begin{cor}
$b^2*X=b*X\implies h^{\star}\le2$ and lengths of
cycles in STD of $\frak A(|X|,X)$ don't exceed 2.
\end{cor}

It's interesting however that for odd $|X|$ the
solutions are the same as for the equation $\frak
b*X=X$. Indeed, for odd $|X|$ the vector $X$ belongs
to the cycle of the baker transformation, i.e.
$X=\frak b^t*X$ for some $t>0$, see
theorem~\ref{th.7}, $\iota(2,|X|)=0$. But because
$\frak b^2*X=\frak b*X$ we can transform $\frak b^t*X$
into $\frak b*X$. However it is not so for even $n$,
in general.

\begin{prop}\label{prop19}
$\frak b*X\neq X\ \&\ \frak b^2*X=\frak B*X\implies
h^{\star}(|X|,X)=2$ or $\hat X$ has an attractor with
length 2.
\end{prop}

{\bf Proof.} In fact, the maximal cycle length is
equal or less than 2. Suppose it is 1. Then if
$\det_2(|X|,X)=1$ we have $X=I$ and hence $\frak
b*X=X$, contradiction. So $\det_2(|X|,X)=0$ and
therefore $h^{\star}(|X|,X)\in\{1,2\}$. However the
case $h^{\star}=1$ leads to the same contradiction in
view of the theorem~\ref{the13}. $\ \ \ \Box$

The rule $[0,0,0,0,0,0,1,1,1,1]$ from the
table~\ref{pic_cyc} is an example: $\frak
b*[0,0,0,0,0,0,1,1,1,1]\neq[0,0,0,0,0,0,1,1,1,1]$ but
$\frak b^2*[0,0,0,0,0,0,1,1,1,1]=\frak
b*[0,0,0,0,0,0,1,1,1,1]$.

The next two easy statements have a clear meaning for
reading of the baker diagram. This is why we place
them here.

\begin{theo}
If $X$ belongs to basin of zero in the baker diagram
for $\frak B^{n}$ then STD for $\frak A(n,X)$ has the
only attractor and it is ${\mathbf 0}$.
\end{theo}

{\bf Proof.} By the condition $\exists i[\frak
b^i*X=\mathbf 0]$. This means that $\exists i\forall
s[\widehat{(\frak b^i*X)}*s=\mathbf 0]$, or every
state $s$ belongs to the basin of zero in STD for $X$.
$\ \ \ \Box$

What are those rules from $\frak B^n$ that belong to
the basin of $\mathbf 0$ of the baker diagram?

The next theorem answers a more general question:

\begin{theo}
$X$ belongs to the basin of a fixed point $Y$ in baker
diagram for $\frak B^n\iff\frak b^{\iota(2,n)}*X=Y$.
\end{theo}

{\bf Proof.} As we know (see, for example,
corollary~\ref{col.6}) $Z=\frak b^{\iota(2,n)}*X$
belongs to a cycle of the baker diagram. And because
the cycle contains only $Y$ we get $Z=Y$.$\ \ \ \Box$

%it was checked - correct
\begin{ex}
{\rm Let $n=6$. Then $\iota(2,n)=1$ and the equation
$\frak b*X=\mathbf 0$ brings $X_0\oplus
X_3=0,X_1\oplus X_4=0,X_2\oplus X_5=0$. From here the
general solution is $X=[a,b,c,a,b,c],a,b,c\in\{0,1\}$.
}
\end{ex}
\ \\

\subsection{Criteria for $h^{\star}=0$ and
$h^{\star}=1$}
\ \\
More complex formulas that include not only $\frak b$
but the binary operation $\boxtimes$, see
definition~(\ref{pre-bak}) have more expressive power.
For example, let's write a criterion of
$h^{\star}(|X|,X)=0$ (or to have non-zero determinant
modulo 2, i.e. ${\det}_2(n,X)=1$).

\begin{theo}~\label{the15}
Let $|X|=n$. Then
\begin{equation}
h^{\star}(n,X)=0\iff\underset{i=0}{\overset{\frak
c(n)-1}{\boxtimes}}\frak b^{\iota(2,n)+i}*X=I.
\end{equation}
\end{theo}

{\bf Proof.} First of all, one can reformulate
$h^{\star}(n,X)=0$ as the proposition that any state
$s$ belong to a cycle of the STD.

Now, if $s$ satisfies the equation
\begin{equation}~\label{ident1}
\hat X^{c^{\star}(n)}*s=s.
\end{equation}
then it belongs to a cycle. Conversely, if any state
$s$ belongs to a cycle of the STD then, taking in
account the theorem~\ref{th10}, we come to the
equation.

Since~(\ref{ident1}) true for any $s$, we get
\begin{equation}~\label{ident2}
\hat X^{c^{\star}(n)}=\hat I.
\end{equation}

What remains is to transform~(\ref{ident2})
equivalently into the equation
\begin{equation*}
\underset{i=0}{\overset{\frak c(n)-1}{\boxtimes}}\frak
b^{\iota(2,n)+i}*X=I
\end{equation*}
i.e. into the form given in the theorem condition.

For that, we need only to note that one can pass from
operator product using the feature of operation
$\boxtimes$ to a $\boxtimes$-composition of vector
notices of the rules. For example, from $\hat X\hat X$
we can pass to $\widehat{X\boxtimes X}$. And then due
to $2^{\iota(2,n)}(2^{\frak c(n)}
-1)=\sum_{i=0}^{\frak c(n)-1}2^{\iota(2,n)+i}$ we can
represent the composition in the form we need by
collection of segments of the composition into blocks
of kind $(\frak b^t*X)$ according to
theorem~\ref{th.5}. $\ \ \ \Box$

\begin{ex}~\label{ex13}
{\rm Let $n=2^m$. Since $\iota(2,2^m)=m,\frak
c(2^m)=1$ the condition of the theorem looks as $\frak
b^m*X=I$. Because $b^m*X$ is swept $\frak
b^m*X=[\sum_{i=0}^mX_i,0,\dots,0]$. Therefore the
solution of the equation is any rule $X$ s.t.
$\sum_{i=0}^mX_i=1$. This means that there are just
$2^{m-1}$ (one half of $\frak B^{n}$) rules with
$h^{\star}=0$. }
\end{ex}

\begin{ex}~\label{ex13a}
{\rm Let $n=2^m3$. We have $\iota(2,2^m)=m,\frak
c(2^m3)=\frak c(3)=2$, and the condition of the
theorem looks as
$$
(\frak b^{m+1}*X)\boxtimes(\frak b^m*X)=I.
$$
 Because $b^m*X$ is $\frak b$-swept, it
 actually has
 3 only non-zero components on positions
 $0\cdot2^m,1\cdot2^m,2\cdot2^{m}$. Let's denote them
 as $x,y,z$ correspondingly. Given with a concrete $m$
 we can easily (using the definition of $\frak b$)
write formulas
 expressing these variables in terms of the components
 of $X$. In particular, when $m=1\ (n=6)$ we get
 \begin{equation}~\label{equa6}
x=X_0\oplus X_3,\ y=X_1\oplus X_4,\ z=X_2\oplus X_5,
 \end{equation}
whereas for $m=2\ (n=12)$ we have
 \begin{eqnarray*}%~\label{equa12}
x&=&X_0\oplus X_3\oplus X_6\oplus X_9,\\
y&=&X_1\oplus X_2\oplus X_4\oplus X_7,\\
z&=&X_5\oplus X_8\oplus X_{10}\oplus X_{11}.
 \end{eqnarray*}
Now, as it was said above,
$b^m*X=[x,0,\dots,0,y,0,\dots,0,z,0,\dots]$ and
$b^{m+1}*X=[x,0,\dots,0,z,0,\dots,0,y,0,\dots]$.
Therefore $(\frak b^{m+1}*X)\boxtimes(\frak
b^m*X)=[x\oplus y\oplus z,0,\dots,0,xz\oplus xy\oplus
yz,0,\dots,0,xz\oplus xy\oplus yz,0,\dots,0]$. Because
the last vector must represent $I$ we come to the
system
 \begin{eqnarray*}
x\oplus y\oplus z&=&1\\
xz\oplus xy\oplus yz&=&0,
 \end{eqnarray*}
or, excluding $x$ from the second equation on the base
of the first,
 \begin{eqnarray}~\label{equa_s}
x&=&1\oplus y\oplus z\\
yz&=&0.
 \end{eqnarray}
What remain is only to replace $x,y,z$ with their
"component meanings" according to~(~\ref{equa6}).
However for that we need to do $n$ certain. We set
$n=6$, and transformed~(\ref{equa_s}) into conditions
for components of $X$ to have $h^{\star}=0$:
 \begin{eqnarray*}%~\label{equa_sc}
1&=&X_0\oplus X_1\oplus  X_2\oplus X_3\oplus X_4
\oplus X_5,\\
0&=&(X_1\oplus X_4)(X_2\oplus X_5).
 \end{eqnarray*}
The numbers of these rules are: 1, 2, 4, 8, 11, 13,
16, 19, 22, 25, 26, 31, 32, 37, 38, 41, 44, 47, 50,
52, 55, 59, 61, 62. Direct computations support this.
}
\end{ex}

The next criterion needs not only equalities by
inequalities too.

\begin{theo}~\label{the14}
Let $|X|=n$. Then $h^{\star}=1$ is true if and only if
 \[ \left\{
 \begin{array}{ccc}
 [\underset{i=0}{\overset{\frak
c(n)-1}{\boxtimes}}\frak
b^{\iota(2,n)+i}*X]\boxtimes X &=& X,\\
 \underset{i=0}
{\overset{\frak c(n)-1}{\boxtimes}}\frak
b^{\iota(2,n)+i}*X &\neq& I.
 \end{array}
\right.
 \]
\end{theo}

{\bf Proof.} First of all we prove that
$h^{\star}(n,X)=1\iff(\det_2(n,X)=0$ and
\begin{equation}~\label{ident3}
\hat X^{c^{\star}(n)}\hat X*s=\hat X*s
\end{equation}
is an identity relatively states $s$.

$\Rightarrow.$ As theorem~\ref{th10} tells us, for any
given $n$ and any $X\in\frak B^n$ lengths of cycles in
STD of $\frak A(n,X)$ must divide number
$c^{\star}(n)=2^{\iota(2,n)}(2^{\frak c(n)} -1)$. Now,
if $h^{\star}(|X|,X)=1$ then the
equation~(\ref{ident3}) is in fact an identity
relatively $s$. Indeed, this is clear not only for the
states $s$ that are not included in any cycle but for
cyclic states too.

$\Leftarrow.$ Let $X$ obeys the
identity~(\ref{ident3}). Then every $X(s)$ belongs to
a cycle of the STD. Therefore the height of the STD
for $X$ can't be bigger 1. And if $\det_2(n,X)=0$ then
the diagram height is not 0.

What remains is to transform equivalently the equation
\begin{equation*}
\hat X^{c^{\star}(n)}\hat X=\hat X
\end{equation*}
into the form given in the theorem condition. This can
be done in the same way as in the previous theorem. $\
\ \ \Box$

\begin{ex}~\label{ex12}
{\rm Let $n=2^m$. Since $\iota(2,2^m)=m,\frak
c(2^m)=1$ the condition of the theorem looks as $\frak
b^m*X\boxtimes X=X$. As we know (see the previous
example also) $\frak
b^m*X=[a,0,\dots,0],a=\sum_{i=0}^mX_i$. The inequality
from the criterion enforces $\det_2(2^m,X)=0$; so
$a=0$. Thus we come to $\mathbf 0\boxtimes X=X$. So
the only rule $X$ for that $h^{\star}(2^m,X)=1$ is
$X=\mathbf 0$. }
\end{ex}

Because of the additional $X$ in the equation of the
last theorem, the calculations become more
complicated, but for $n=6$ it is quite doable even by
hands. In this way we found the list of all rules of
length 6 that have $h^{\star}=1$: 0, 5, 10, 15, 17,
20, 21, 30, 34, 39, 40, 42, 51, 57, 60. So the rules
that do not occur in this list and the list of the
example~\ref{ex13a} have the height 2, since according
to our upper estimation 2 is the upper limit
for $n=6$. \\

Of course, this line of criteria can be continued,
i.e. one can formulate analogously criteria for
$h^{\star}$ to be equal to given a number $k$.
However, the complexity would increase and therefore
the computational aspect of these expressions deserves
a discussion.

Coming to the computational aspect of these results in
general, let's estimate the complexity of the
computation setting by formula
$\underset{i=0}{\overset{\frak
c(n)-1}{\boxtimes}}\frak b^{\iota(2,n)+i}*X$. Since we
deal with boolean strings, $\frak b*X$ can be
calculated for $\mathcal O(n)$ time, where $n=|X|$.
Therefore in sum to compute all operands of
$\underset{i=0}{\overset{\frak c(n)-1}{\boxtimes}}$ we
need no more than $\mathcal O(nc^{\star}(n))$ of time.
Also every $\boxtimes$ with $n$-long boolean strings
needs no more than $\mathcal O(n^2)$ time. Therefore
we estimate the general time expenses as $\mathcal
O(\frak c(n)\cdot n^2)$. This is comparable with the
time one needs to calculate the rank of a $n\times
n$-matrix with boolean elements. According to [1] the
average value $\overline{\frak c(n)}$ of $\frak c(n)$
grows as $\mathbf o(n)$, i.e.
\begin{equation}~\label{equa1}
\overline{\lim_{n\to\infty}}\frac{1}{n^2}\sum_{1\le
i\le n}\frak c(n)=0.
\end{equation}
This mean that the method, suggested in
theorem~\ref{the15}, to calculate our determinants
modulo 2 is more effective in average than the method
using rank computation.

Anyway, in case of small $\frak c(n)$ the suggested
method to compute the determinant modulo 2 has the
good practical efficiency.

Therefore, the results presented by
theorems~\ref{the15},\ref{the14} a certain theoretic
and computational value.\\

\subsection{Description of the
basin of $\mathbf 0$ in STD}
\ \\
It is well known with given rule $X$, how to write the
system of linear equations s.t. its solutions are
pre-images of $\mathbf 0$. However to describe the
whole basin of $\{\mathbf 0\}$ we need to solve
equation systems of kind $X^i*s=\mathbf
0,i=1,2,\dots$. The next result suggests a single
equation, describing $\mathbf 0$-basin of $X$. We
recall that $\circlearrowleft$ is reversion of
sequence, and $\sigma$ is the cyclic shift to right on
one position.

\begin{theo}
$s$ belongs to the basin of $\{\mathbf 0\}$ in STD for
$\frak A(n,X)\iff s\boxtimes(\sigma*(\frak
b^{\iota(2,n)}*X)^{\circlearrowleft})=\mathbf 0$.
\end{theo}

{\bf Proof.} First of all, as we know,
$h^{\star}(n,X)\le2^{\iota(2,n)}$. Hence if $Z=\frak
b^{\iota(2,n)}*X$ then ($s$ belongs to the basin of
zero$\iff\hat Z*s=\mathbf 0$).

Now, all what remain to do is to write the circulant
matrix equation $C(Z)*s=\mathbf 0$ in terms of our
operations $\boxtimes,\circlearrowleft,\sigma$. For
that we pass firstly to $s*C^T(Z)$, where $C^T$ is the
transposed circulant $C$, and then apply $\boxtimes$
instead of the right multiplication of the vector on
the matrix. So $\sigma*(Z^{\circlearrowleft})$ present
the first column (row) of the circulant $C(Z)\
(C^T(Z))$. $\ \ \ \Box$\\

\subsection{Determinant reduction}
\ \\
For the next result it is convenient to introduce the
operation of $\frak b$-compression $X/{\frak b}$ of
given sequence $X$ as following. First of all we pass
from $X=[x_0,\dots,x_{n-1}]$ to $\frak
b^{\iota(2,n)}*X=[z_0,0,\dots,z_1,0,
\dots,z_{\frac{n}{\iota(2,n)}},0\dots]$. Here 0 occupy
positions $j$ s.t. $2^{\iota(2,n)}\nmid j$. The
positions $j$ which are multiple of $2^{\iota(2,n)}$
are occupied by $z_i$. At last, $X/{\frak
b}=_{df}[z_0,z_1,\dots,z_{\frac{n}{\iota(2,n)}}]$.

\begin{theo}\label{th12}
$\det_2(n,X)=\det_2(\frac{n}{\iota(2,n)},X/{\frak
b})$, i.e. in other words the determinant modulo 2 of
the any rule $X$ coincides with the determinant modulo
2 of the result $\frak b$-compression applied to $X$.
\end{theo}

{\bf Proof.} We suggest two proofs.

(I) Let's start from the theorem~\ref{the15}. It can
be reformulated as ($n=|X|$):
\begin{equation*}
{\det}_2(n,X)=1\iff\underset{i=0}{\overset{\frak
c(n)-1}{\boxtimes}}\frak b^{\iota(2,n)+i}*X=I.
\end{equation*}
Yet, $\frak b^{\iota(2,n)+i}*X=\frak b^i*Y$ if $Y$
denotes $\frak b^{\iota(2,n)}*X$. On the other side,
as we know, all components of $Y$ which numbers are
not divisible by $2^{\iota(2,n)}$ are equal to 0. The
same is true for $\frak b^i*Y$. And according to
lemma~\ref{3a} all components with numbers not
divisible by $2^{\iota(2,n)}$ of the left side
\begin{equation}~\label{eq_det}
\underset{i=0}{\overset{\frak c(n)-1}{\boxtimes}}\frak
b^i*Y
\end{equation}
of the equation in the equivalence are 0.

The following discourse in essence repeats the
reduction lemma~\ref{lem6}. Let's pass from $Y$ to a
rule $y$ by cancelling all components of $Y$ that have
numbers being not divisible by $\iota(2,n)$. Clearly,
$y$ is $\frak b$-compression of $X$. Sure, $y\in\frak
B^{n'}$ where $n'=\frac{n}{\iota(2,n)}$. Due to
lemma~\ref{3a} if $y$ is $\frak b$-compression of $X$
then $\frak b*Y$ also can be transformed into $y'$
that is $\frak b$-compression of $\frak b*X$, and
$y'=\frak b*y$. Therefore the equality
$$
\underset{i=0}{\overset{\frak c(n)-1}{\boxtimes}}\frak
b^{\iota(2,n)+i}*X=I
$$
is true if and only if
$$
\underset{i=0}{\overset{\frak c(n)-1}{\boxtimes}}\frak
b^i*y=I
$$
is true. Of course, last $I$ denote the sequence
$[1,0,\dots]\in\frak B^{n'}$.

By theorem~\ref{the15} the last equality is equivalent
$\det_2(n',y)=1$. So,
$\det_2(n,X)=1\iff\det_2(n',y)=1$ if only $y$ is
$\frak b$-compression of $X$.

(II) The second proof is based on the classical
formula for the determinant of a circulant matrix.

As it's well known determinants of linear operator $L$
and every its degree $L^i,i\ge1,$ are equal or not
equal to 0 simultaneously. For the boolean field and
determinants modulo 2 this leads to the possibility to
replace equality to 0 with equality. Therefore
$\det_2(|X|,X)=\det_2(|X|,\frak b^m*X)$. If we set
$m=\iota(2,|X|)$ then we get $\frak b$-compression
$X/{\frak b}$ of $X$. Now we write out the determinant
$\delta$ of the circulant matrix, produced by $\frak
b^{m}=[z_0,z-1,\dots,z_{|X|}]$, using the formula
(14.312) from [3, p.1068]:
\begin{equation}
\delta=\prod^{|X|}_{j=1}\sum_{i=0}^{|X|-1}z_iw^i_j.
\end{equation}
Because, $z_i=0$ if $2^m\nmid i$ we can write
$\delta=\prod^{|X|}_{j=1}
\sum_{i=0}^{|X|/2^m-1}z_{i2^m}w^{i2^m}_j$. When $i$
runs the list $0,\dots,|X|/m-1$ the number
$w^{i2^m}_j$ runs subsequent roots
$e^{\frac{2\pi}{i}},i=0,\dots,|X|/2^m-1$ of degree
$|X/{\frak b}|$ of 1. So we come to
$\delta=\delta_1^{2^m}$ where $\delta_1$ is the
determinant of $\frak b$-compression $X/{\frak b}$.
Now the conclusion is obvious. $\ \ \ \Box$

This theorem can help in case when the compression
factor $\frac{n}{\iota(2,n)}$ is sufficiently big
comparing with $n$. In particular, as we already know,
$\det_2(2^k,X)=X_0\oplus X_1\oplus\dots\oplus
X_{2^k-1}$.

\section{References}

1. O.Martin, A.Odlyzko, S.Wolfram. {\em Algebraic
Properties of Cellular Automata}, in the book {\em
Cellular Automata and Complexity} by S.Wolfram, 1994,
pp. 71-113.

2. B.Voorhees, C.Beauchemin. {\em Point Mutations and
Transitions Between Cellular Automata Attractor
Basins.} Preprint, {\tt arXiv:nlin.CG/0306033} v1, 17Jun
2003.

\end{document}